\documentclass[12pt]{article}
\usepackage{amssymb}
\usepackage[pdftex]{graphicx}

\title{The straight line, the catenary, the brachistochrone, the circle, and Fermat}
\author{Raul Rojas\\Freie Universit\"at Berlin}
\date{January 2014}  
\begin{document}
\maketitle

\begin{abstract}
This paper shows that the well-known curve optimization problems which lead to the straight line, the
catenary curve, the brachistochrone, and the circle, can all be
handled using a unified formalism. Furthermore, from the general differential equation fulfilled
by these geodesics, we can guess additional functions and the required metric. The parabola, for example,
is a geodesic under a  metric guessed in this way.
Numerical solutions are found for the curves corresponding to geodesics in the 
various metrics using a ray-tracing approach based on Fermat's principle.
\end{abstract}

\section{Motivation}

Students of mathematics, physics, or engineering, get acquainted with problems of variational calculus
during their first years in college. In my case, I learned to find the catenary curve in a mechanics course.
We learned about the brachistochrone in a further course about theoretical mechanics (where the Euler-Lagrange equation
plays a major role). The properties of the circle were studied in a geometry class, and I learned to use
semicircles as models for the lines in hyperbolic geometry after reading a book on non-Euclidean geometries.
Since the analytic solutions for each variational problem  look very different (the catenary is a sum of exponentials,
while the brachistochrone can be expressed parametrically using sines and cosines), and the procedure used to find the solution
also changes significantly from book to book, it is not immediately obvious that we are esentially dealing  with the same
problem, even when some of them appear in a book solved one after the other \cite{fox}. That is what will be shown here. First we look at  each optimization problem. We will reduce
them  to a unified formulation, and we will then solve them analytically and numerically.

\begin{figure}[htb]
\centering
\includegraphics[width=\linewidth]{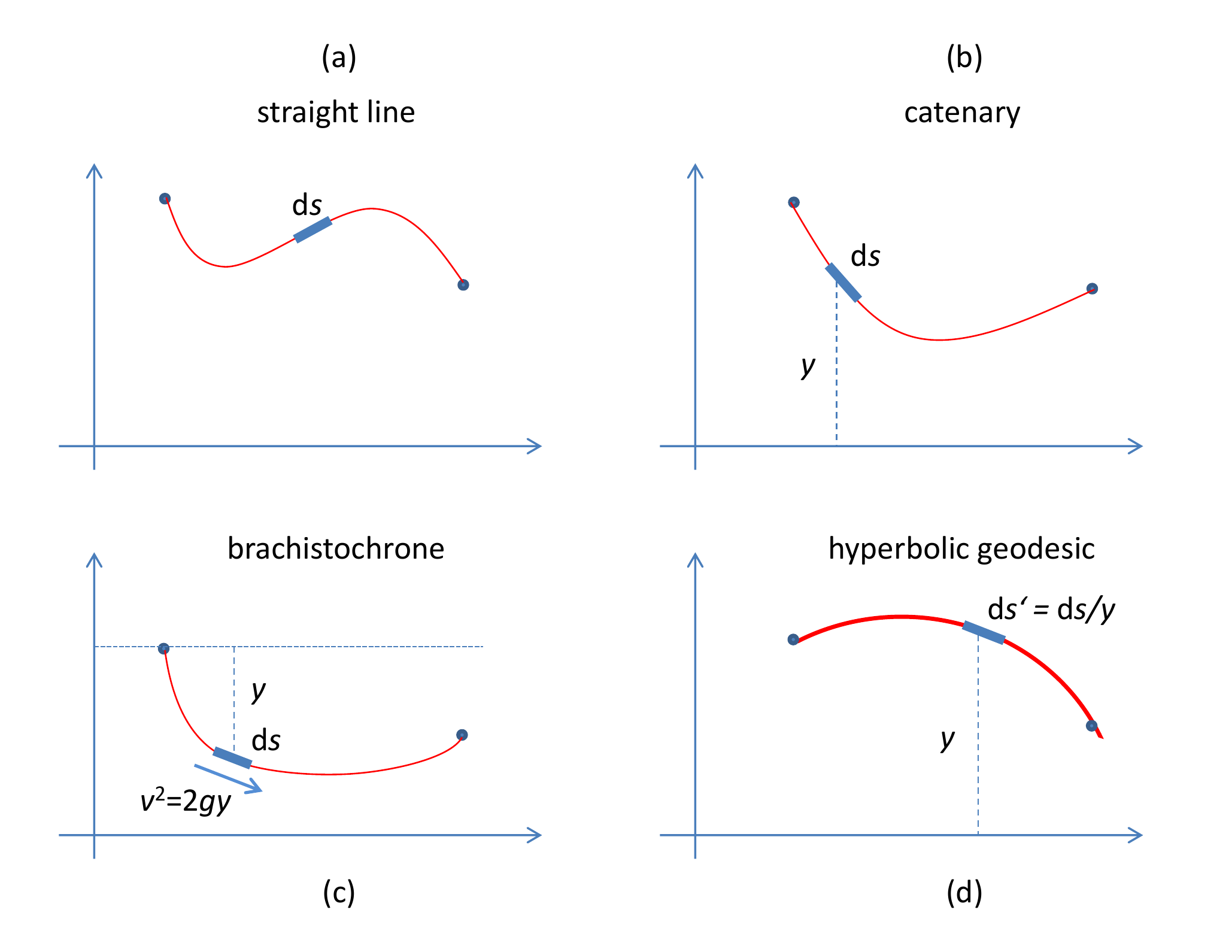}
\caption{\label{Abb1} The four optimization problems: (a) finding the shortest path in Euclidean geometry, (b) the shape of a hanging chain, (c) the path of fastest descent, and (d) the shortest path in hyperbolic geometry.}
\end{figure}

\noindent{\bf Straight line}

Our first example is finding the shortest curve $f$ between two points (in Euclidean geometry, Fig.~\ref{Abb1}a). The differential arc length is given by ${\rm d}s^2={\rm d}x^2+{\rm d}y^2$, or written in another way, ${\rm d}s=\sqrt{1+y'^2}{\rm d}x$. The total length $L$ between two points $(x_1,y_1)$ and $(x_2,y_2)$ along a curve $f$ in the Euclidean plane is given by
$$L = \int_{x_1}^{x_2}  {\rm d}s \hspace{1cm}{\rm subject\ to}\hspace{1cm} f(x_1)=y_1, f(x_2)=y_2.$$

\noindent{\bf Catenary}

In the case of a chain hanging from two given points, what we want to minimize is the total potential energy
of the chain (Fig.~\ref{Abb1}b). A differential piece of the chain, of length ${\rm d}s$ has mass ${\rm d}m=\rho {\rm d}s$, where $\rho$ is the mass density of the chain. The height $y$ of each differential segment, multiplied by its mass, and the gravitational constant $g$, is equal to the potential energy of the segment. The total potential energy $E$ of a chain hanging with the shape $f$ is given by

  $$E = \int_{x_1}^{x_2} y g(\rho {{\rm d}}s)= g\rho\int_{x_1}^{x_2} y {\rm d}s\hspace{1cm} {\rm subject\ to}\hspace{0.5cm} f(x_1)=y_1, f(x_2)=y_2.$$

\noindent{\bf Brachistochrone}

The Brachistochrone is the curve $f$ for a ramp along which an object can slide from rest at a point $(x_1,y_1)$ to a point $(x_2,y_2)$ in  minimal time (Fig.~\ref{Abb1}c). Since the speed of the sliding object is equal to $\sqrt{2gy}$, where $y$ is measured vertically downwards from the release point, the differential time it takes the object to traverse the arc ${\rm d}s$ at that speed is ${\rm d}s/\sqrt{2gy}$. The total travel time $T$ is  given by
$$T = \int_{x_1}^{x_2}{(2gy)^{-1/2} \rm{d}s}= (2g)^{-1/2}\int_{x_1}^{x_2}{y^{-1/2} \rm{d}s} \hspace{0.5cm} {\rm subject\ to} \hspace{0.5cm} f(x_1)=y_1, f(x_2)=y_2$$.

\noindent{\bf Hyperbolic geodesics}

A model for hyperbolic geometry is a half-plane in which the length of differential segments varies
vertically  inversely proportional to $y$ (Fig.~\ref{Abb1}d). That means: in the hyperbolic metric for this halfplane a segment
has length $ds'=ds/y$, where $ds$ is the usual Euclidian metric. In that case finding the shortest path between
two points has to be made considering the new metric. The length of the path L between two points along the curve $f$ is given by
 $$L=\int_{x_1}^{x_2} y^{-1} {\rm d}s  \hspace{1cm} {\rm subject\ to} \hspace{1cm}f(x_1)=y_1, f(x_2)=y_2.$$

\section{The general formulation}

From what has been said it is obvious that the general formulation of the four problems above has the following form: Find the curve $y$ that minimizes
\begin{equation}
k\int_{x_1}^{x_2} {y^{\alpha} {\rm d}s}=k\int_{x_1}^{x_2}y^{\alpha}{\sqrt{1+y'^2}}{\rm d}x \hspace{1cm} {\rm subject\ to} \hspace{1cm} f(x_1)=y_1, f(x_2)=y_2,
\label{general}
\end{equation}
for $\alpha=0$ (Euclidean geometry),$\alpha=1$ (catenary),$\alpha=-1/2$ (brachistochrone), and $\alpha=-1$ (hyperbolic geodesic), and where $k$ is a constant.

The general method for finding a solution to this problem of variational calculus would be to use the Euler-Lagrange equation \cite{lanczos}. Given the problem of finding an optimal value for an integral of the form
$$\int_a^b L(x,y,y')dx$$
we can solve the differential equation
\begin{equation}
\frac{\partial L }{\partial y} - \frac{\rm d}{{\rm d}x} \frac{\partial L }{\partial y'} = 0
\label{euler}
\end{equation}
in order to find a solution.

Since in our case the argument of the integral operator $L=y^{\alpha}\sqrt{1+y'^2}$ has  partial derivative relative to $x$ equal to zero, a simplification of the Euler-Lagrange equation can be used, the so-called Beltrami's identity:
\begin{equation}
L - y'\frac{\partial L }{\partial y'} = C
\label{beltrami}
\end{equation}
Since
$$\frac{\partial }{\partial y'}\sqrt{1+y'^2} = \frac{y'}{\sqrt{1+y'^2}}$$
applying Beltrami's identity (Eq.~\ref{beltrami}) to our general formulation in Eq.~\ref{general} leads to the following differential equation
$$y^{\alpha}{\sqrt{1+y'^2}} - \frac{y^\alpha y'^2}{\sqrt{1+y'^2}} = C$$
which can be simplified to
\begin{equation}\frac{y^\alpha}{\sqrt{1+y'^2}}=C
\label{gen1}\end{equation}
or equivalently
\begin{equation}
\label{eqdif}
\frac{y^{2\alpha}}{{1+y'^2}}=\frac{y^{2\alpha}{\rm d}x^2}{{{\rm d}x^2+{\rm d}y^2}}=C^2,
\end{equation}
which is a first order nonlinear ordinary differential equation.

\noindent {\bf Proof of Beltrami's Identity}

From the Euler-Lagrange equation we know that
\begin{equation}
\label{step}
\frac{\partial L}{\partial y} = \frac{{\rm d}}{{\rm d}x}\frac{\partial L}{\partial y'}
\end{equation}
Using the chain rule we can express the derivative of $L$ according to $x$ as follows:
$$\frac{{\rm d} L}{{\rm d} x} = y'\frac{\partial L}{\partial y}+y''\frac{\partial L}{\partial y'}+\frac{\partial L}{\partial x}$$
In the case that ${\partial L}/{\partial x}=0$, we can simply write
$$\frac{{\rm d} L}{{\rm d} x} - y'\frac{\partial L}{\partial y}-y''\frac{\partial L}{\partial y'}=0.$$
Substituting the equivalent of ${\partial L}/{\partial y}$ from Eq.~\ref{step} in the second term above and grouping, we obtain
$$\frac{{\rm d} L}{{\rm d} x} - (y'\frac{{\rm d}}{{\rm d}x}\frac{\partial L}{\partial y'}+y''\frac{\partial L}{\partial y'})=0$$
But then the last two terms can be rewritten (using the differentiation product rule) as
$$\frac{{\rm d} L}{{\rm d} x} - \frac{{\rm d}}{{\rm d}x}\left(y'\frac{\partial L}{\partial y'}\right)=\frac{{\rm d}}{{\rm d}x}\left(L- y'\frac{\partial L}{\partial y'}\right)=0$$
and since the derivative is zero, we conclude that
$$L -y'\frac{\partial L}{\partial y'}=C$$
for a certain constant C.

It is easy to check that the solutions for the catenary, brachistochrone and circle fulfill the differential
equations derived above (Eq.~\ref{gen1} and Eq.~\ref{eqdif}).

{\centering 
\begin{tabular}{|l| c| l |l | c|}
\hline 
Problem &$\alpha$& Parametric solution & ${\rm d}{x}$, ${\rm d}{y}$ & $C^2$\\
\hline 
Line & 0 &$\begin{array} {lll} x & = t&   \\  y & = at + b&  \end{array}$ & $\begin{array} {lll} {\rm d}x & = 1&   \\  {\rm d}y & = a&  \end{array}$ & $\frac{1}{1+a^2}$\\ 
\hline 
Catenary & 1 &$\begin{array} {lll} x & = t&   \\  y & = {\rm cosh} (t)&  \end{array}$ & $\begin{array} {lll} {\rm d}x & = 1&   \\  {\rm d}y & = {\rm sinh} (t)&  \end{array}$ & 1\\ 
\hline 
Brachistochrone & $-1/2$ &$\begin{array} {lll} x & = r(t-{\rm sin} (t))&   \\  y & = r(1-{\rm cos} (t))&  \end{array}$ & $\begin{array} {lll} {\rm d}x & = r(1-{\rm cos} (t))&   \\  {\rm d}y & = r {\rm sin} (t)&  \end{array}$ & $\frac{1}{2r}$\\ 
\hline 
Hyperbolic & $-1$ &$\begin{array} {lll} x & = r{\rm cos}(t)&   \\  y & = r{\rm sin}(t)&  \end{array}$ & $\begin{array} {lll} {\rm d}{x} & = -r{\rm sin}(t)&   \\  {\rm d}{y} & =r{\rm cos}(t) &  \end{array}$ & 1\\ 
\hline 
\end{tabular}}

The interesting thing is that having the differential equation Eq.~\ref{eqdif} we can now guess solutions, and from them infer the corresponding metric. The parabola $y=\frac{1}{4}x^2+1$, for example, fulfills the differential equation for $\alpha=1/2$, as the reader can verify. We can then say that for a metric of the form $\sqrt{y}{\rm d}s$ the parabola is the appropriate geodesic (the curve of minimal length).

\section{Numerical solutions}

The variational problem can be solved numerically or algebraically. One way to go is to use an ``analog computer'' designed 
for minimizing soap bubbles, the other is by using an equation solver. 

\begin{figure}[htb]
\centering
\includegraphics[width=0.6\linewidth]{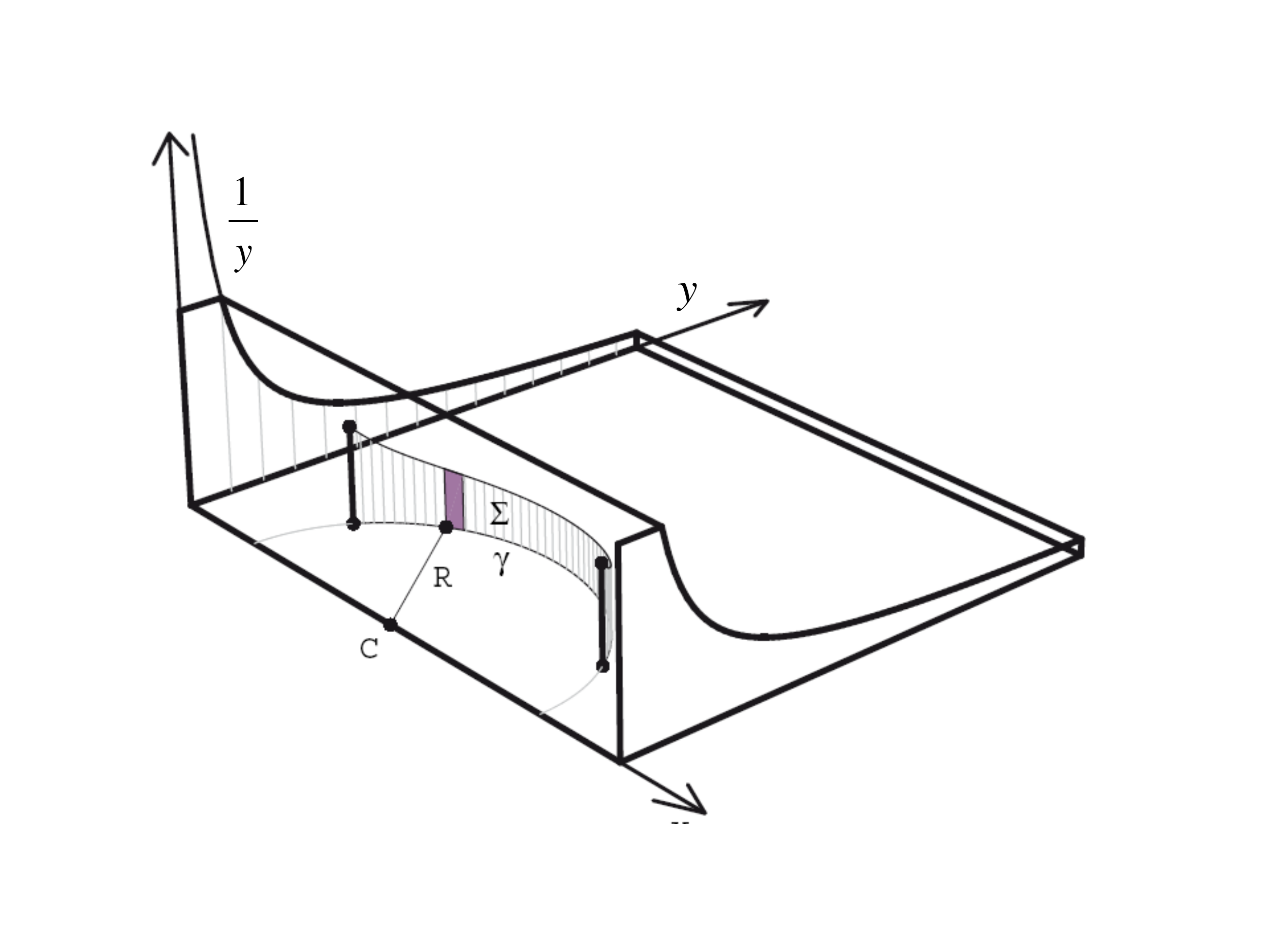}
\caption{\label{soap}A soap-film computer (from \cite{criado}). The roof of the box follows the function $1/y$. The coordinate axis for $y$
is shown in the diagram. Two vertical pins enclose a soap film $\Sigma$. The film of minimal area has the shape of an arc of a circle of radius $R$}
\end{figure}

Fig.~\ref{soap} shows the general approach applied by Criado and Alamo \cite{criado}. A box is built and the top of the box is curved in such a way that the height of this ``roof'' follows a certain function ($1/y$ in the case of Fig.~\ref{soap}). If the box is submerged in soap water and is pulled out, a soap film can build between the two vertical pins, and the bottom and top of the box. The area of the soap film $\Sigma$ is proportional to the length of the film, and its height at every point along the film path. It is easy to see that this corresponds to  a geodesic problem. For the metric $1/y$ the geodesic is an arc of a circle. The area of the film $\Sigma$ represents the integral of $y^{-1}ds$, because the area of a differential portion of the soap film is precisely $y^{-1}ds$. Therefore the soap film of minimal area will follow an arc of a circle, as can be confirmed experimentally. Any other metric can be ``processed'' by changing the shape of the top of the box, and the corresponding geodesic can be found.

Alternatively, the Mathematica solver (Wolfram Alpha)  for differential equations can be used in order to obtain solutions for Eq.~\ref{eqdif}. The analytic solutions for $\alpha=-1,-1/2,0,1/2,1$ are readily obtained. It is also possible to obtain some solutions for some special values of $\alpha$, such as $\alpha=0.75/2$, in which case the following solution, expressed as an implicit function of $x$ and $y$, is obtained:
$$c_1+0.5 x = 4 y^{\frac{1}{4}}+\frac{4}{3} {\rm log}(1-y^{\frac{1}{4}})-\frac{2}{3} {\rm log}(\sqrt{y}+y^{\frac{1}{4}}+1)-\frac{1}{\sqrt{3}}4 {\rm tan}^{-1}\left(\frac{2 y^\frac{1}{4}+1}{\sqrt{3}}\right).$$

It is interesting that the ``historical'' cases, for $\alpha=-1,-\frac{1}{2},0,\frac{1}{2}, 1$, are precisely those which have a nice closed solution. Other similarly simple functions could be possibly found for some other values of $\alpha$ using the differential equations solver.

\section{Fermat's Principle}

Johann Bernoulli solved the brachistochrone problem in 1696/97 casting it as the problem of computing the travel time of light crossing through a medium where its speed is changing \cite{sussmann}. The refraction of light passing from a medium A to a medium B can be 
computed in optics applying Fermat's principle. This states that the travel time from a point A to a point B should be minimal, and if the velocity of light is different in medium A and medium B, then the optimal travel direction is such that
$$\frac{sin(\theta_1)}{v_1}=\frac{sin(\theta_2)}{v_2}$$
where $\theta_1$ and $\theta_2$ represent the angles with respect to the normal through the boundary between the two media and $v_1$, $v_2$ are the respective velocities of light in each medium. This means that Fermat's principle can be reduced to the statement that ${sin(\theta)}/{v}$ 
is a constant when a light ray is going through different media (and $\theta$ and $v$ have the interpretations given above).

\begin{figure}[htb]
\centering
\includegraphics[width=0.6\linewidth]{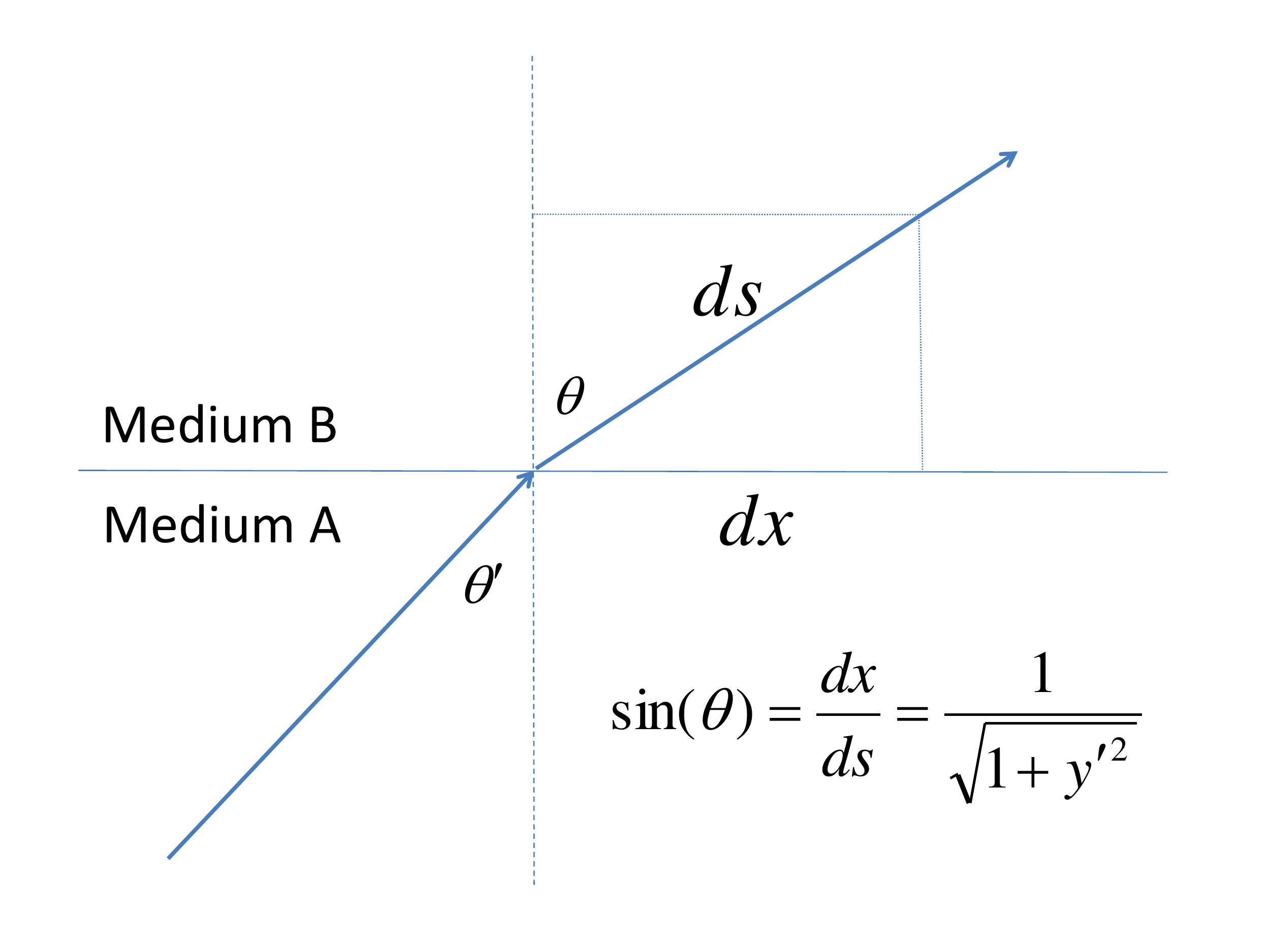}
\caption{\label{fermat}Fermat's principle. A light ray crosses from medium A to medium B ansd is refracted.}
\end{figure}

The sine of the angle with the vertical can be computed as shown in Fig.~\ref{fermat}. It is equal to
$$sin(\theta)=\frac{dx}{ds}=\frac{dx}{\sqrt{dx^2+dy^2}}=\frac{1}{\sqrt{1+y'^2}}.$$

Now, if we go back to the brachistochrone problem, the gerneral formulation for the differential equation to be solved, obtained from Beltrami's identity (Eq.~\ref{gen1}), was
$$\frac{y^{\alpha}}{\sqrt{1+y'^2}}=C$$
for $\alpha=-1/2$. Since the velocity $v$ of the falling object is proportional to $\sqrt{y}$, this means that the expression obtained from Beltrami's identity is nothing else but Fermat's principle! This is so because $v\sim\sqrt{y}$ and
$$\frac{y^{-1/2}}{\sqrt{1+y'^2}}= \frac{sin(\theta)}{\sqrt{y}}=C$$

\begin{figure}[htb]
\centering
\includegraphics[width=0.6\linewidth]{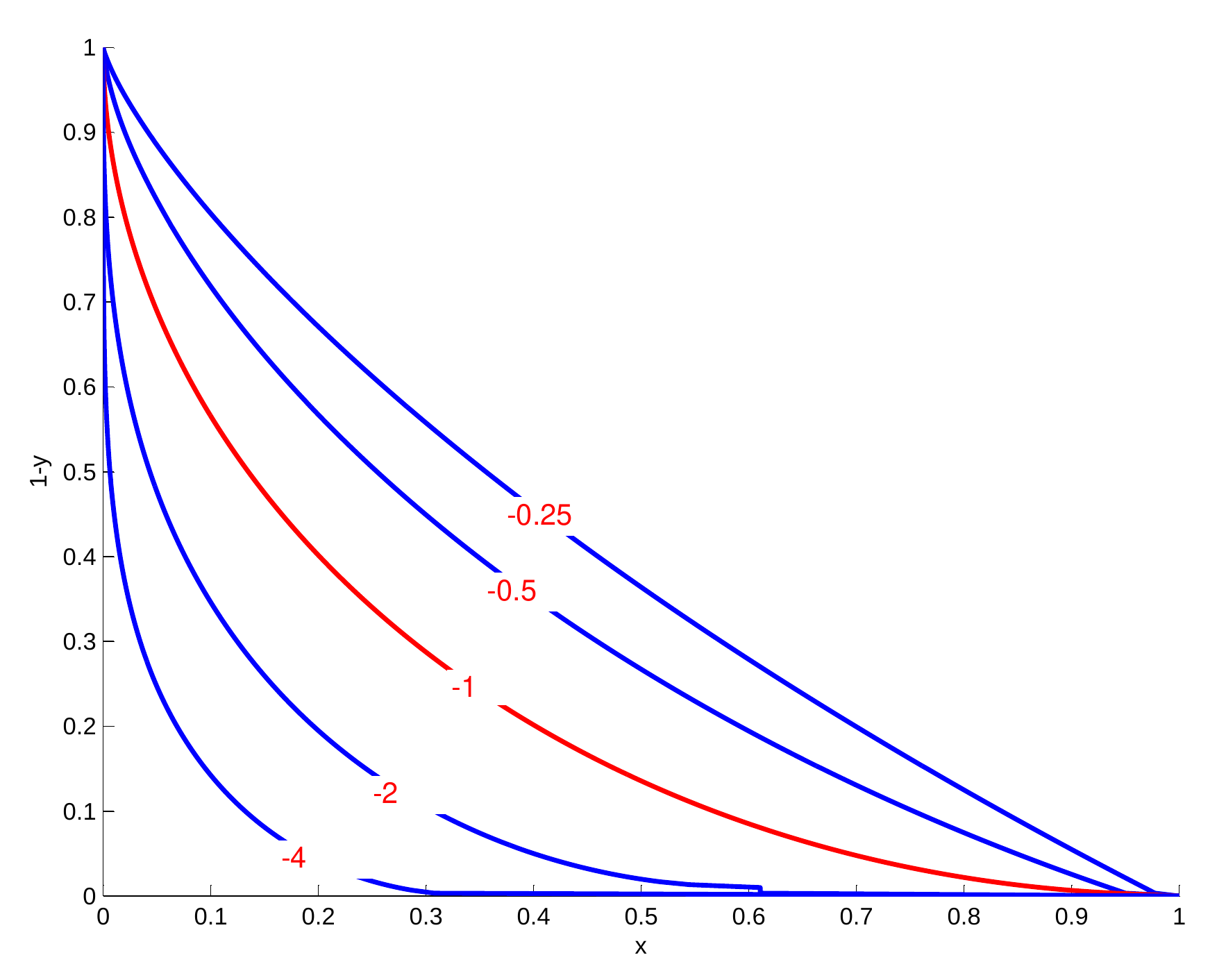}
\caption{\label{geodesics}Geodesics for $\alpha=-0.25,-0.5,-1,-2,-4$ (from top to bottom curve). The circle ($\alpha=-1$) has been plotted red. The curve immediately above the circle is the brachistochrone. The vertical axis is $1-y$, so that the brachistochrone can be compared with the circle, the hyperbolic geodesic.}
\end{figure}

Generalizing, we now interpret Eq.~\ref{gen1} in this form
$$\frac{y^{\alpha}}{\sqrt{1+y'^2}}= \frac{sin(\theta)}{y^{-\alpha}}=C,$$
This corresponds to the general case in which the ``velocity'' of the particle in each medium is proportional to $y^{-\alpha}$, that is, proportional to the generalized metric. It is as if we had modified the gravity law to fit a  new metric. The corresponding geodesics can be obtained solving the differential equation as before, or just numerically, applying Fermat's principle to the graph of the solution of the differential equation. The numerical solution can then be found in a straightforward manner, using ray tracing for the generalized metric. Fig.~\ref{geodesics} shows the results for various negative values of alpha (the vertical axis is $1-y$ in order to make the brachistochrone comparable with the circle, the hyperbolic geodesic).

\section{Surfaces of revolution and designer geodesics}

If we now let the the catenary rotate around the horizontal axis, we obtain a kind of curved cylinder. The surface of the side of the cylinder is minimal, since the curve that minimizes
  $$E = g\rho\int_{x_1}^{x_2} y {\rm d}s\hspace{1cm} {\rm subject\ to}\hspace{0.5cm} f(x_1)=y_1, f(x_2)=y_2$$
also minimizes
$$E = 2\pi\int_{x_1}^{x_2} y {\rm d}s\hspace{1cm} {\rm subject\ to}\hspace{0.5cm} f(x_1)=y_1, f(x_2)=y_2.$$
This last integral represents the surface of a curved-cylinder of revolution with radius $y_1$ and $y_2$ at the basis and top. The minimal surface of revolution of this kind is called the catenoid, because it is obtained rotating a catenary curve.

\begin{figure}[htb]
\centering
\includegraphics[width=0.6\linewidth]{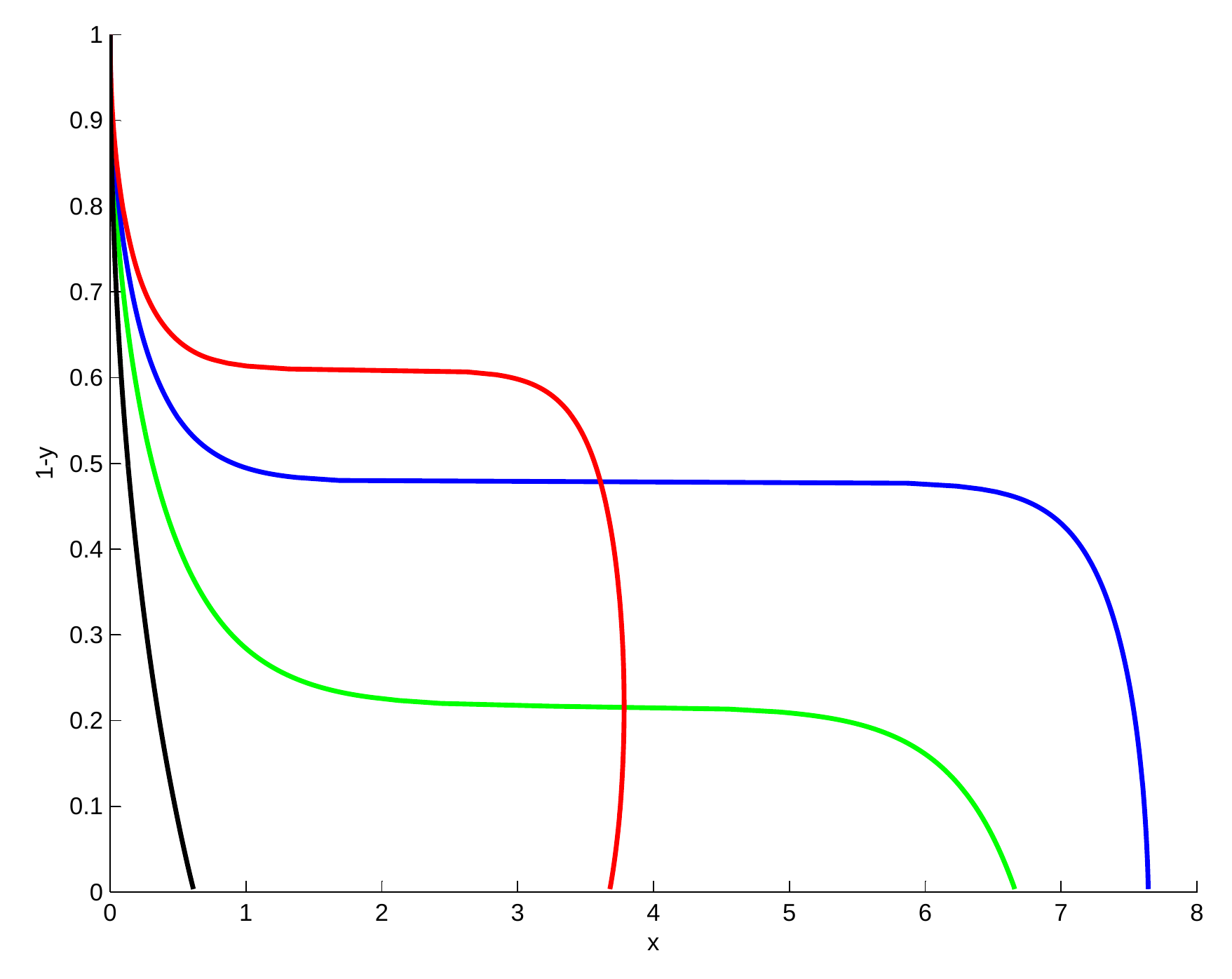}
\caption{\label{designer}Geodesics for metrics of the form $1/\sin(a(1-y))$, for $a=1$ (black), $a=2$ (green), $a=3$ (blue), and $a=4$ (red). The red line curves back because the index of refraction achieves negative values near to the axis, for the metric selected. No boundary conditions were imposed on the solutions.}
\end{figure}

It is also possible to force almost any monotonous curve to be a geodesic. Fig.~\ref{designer} shows an example where the metric has been chosen not proportional to $(1-y)^{\alpha}$ but to functions of the form $1/\sin(a(1-y))$, so that the velocity of ``light'' through such a medium is forced to be proportional to $C\sin{(a(1-y))}$. In Fig.~\ref{designer} the values of $a$ used are $1,2,3$ and $4$.

We can also just invert Fermat's principle: given a desired geodesic, determined by a sequence of angles $\theta_1,\theta_2,\ldots,\theta_n$ relative to the vertical, of segments $ds$, anchored at heights $y_1,y_2,\ldots,y_n$, Fermat's principle tells us that we need velocities in each segment, such that $v_i=\sin(\theta_i)/C$, for $i=1,\ldots,n$.
The necessary refraction indices for each layer can be computed from the relative change in angle $\theta_i$ going through each layer.

In some cases we can even find the metric for a desired geodesic function analytically. Assume that we want the hyperbola to be a geodesic for some metric $g(y)$, which does not depend explicitly on $x$ (so that we can apply Beltrami's identity).
In this case we want to minimize 
$$
\int_{x_1}^{x_2} g(y){\sqrt{1+y'^2}}{\rm d}x
$$
disregarding the boundary conditions. Substituting the metric $g(y)$ for the metric $y^\alpha$ in the differential Eq.~\ref{gen1} we are left with the condition
\begin{equation}\frac{g(y)}{\sqrt{1+y'^2}}=C.
\label{hyper}\end{equation}
Since for the hyperbola $y=1/x$ the derivative is $y'=x^{-2}=y^2$, but then ${\sqrt{1+y'^2}}=\sqrt{1+y^4}$, we only need to set $g(y)={\sqrt{1+y^4}}$ so that Eq.~\ref{hyper} is fulfilled for $C=1$.

In general, any function $y$ whose derivative $y'$ can be written in the form $y'=h(y)$, where $h$ is a function of $y$, but not explicitly of $x$, can be a geodesic for the metric $g(y)=C\sqrt{1+h(y)^2}$, where $C$ is a constant. This is possible for $y=e^x$, or $y=\sin^2{x}$, for example. If boundary conditions are considered, the problem becomes harder, and a solution could not exist, within the family of functions selected at the beginning.

\section{Final comments}

All the variational problems discussed here modify the spatial metric in the plane using a function of the height $y$ and independent of $x$. This makes possible the application of Beltrami's identity to all these problems.
The metrics discussed are of the type $y^\alpha {\rm d}s$ or $g(y){\rm d}s$. Any value of $\alpha$ can be conceivably used, as we did when solving numerically. The physical interpretation for such curious metrics correspond to a plane with a continuously vertically varying index of refraction (positive or negative). Since the metric does not depend on $x$, the plane has been filled with horizontal layers with varying index of refraction (like a  stack along the $y$-direction of different transparent materials). This makes light rays curve along the geodesics of the resulting metric.

Materials with a varying index of refraction have been discussed in optics and have become fashionable in the context of ``cloaking devices'' that curve light rays around a disguised object \cite{sarbort}. Such materials represent a form of ``analog computer'' for geodesics ... without the soap water.

\end{document}